\documentclass[11pt]{article}
\usepackage{graphicx}

\newfont{\bb}{msbm10 at 10pt}
\def\r{\hbox{\bb R}}

\newcommand{\e}{\hbox{\bf E}}
\newcommand{\h}{\hbox{\bf H}}
\newtheorem{theorem}{Theorem}[section]
\newtheorem{corollary}{Corollary}
\newtheorem{definition}{Definition}[section]
\newtheorem{remark}{Remark}

\begin{document}

\title{Spacelike surfaces in Minkowski space satisfying  a linear relation between their principal curvatures}
\author{\"{O}zg\"{u}r Boyac\i o\u{g}lu Kalkan \\
Mathematics Department\\
Afyon Kocatepe University\\
Afyon 03200 Turkey\\
email: bozgur@aku.edu.tr\\
 \vspace*{.5cm}\\
Rafael L\'opez\thanks{%
Partially supported by MEC-FEDER grant no. MTM2007-61775 and Junta de Andaluc\'{\i}a grant no. P06-FQM-01642.}\\
Departamento de Geometr\'{\i}a y Topolog\'{\i}a\\
Universidad de Granada\\
18071 Granada, Spain\\
email: rcamino@ugr.es}
\date{}
\maketitle

\begin{abstract}
In this work, we consider spacelike surfaces in Minkowski space $\hbox{\bf E}%
_{1}^{3}$  that satisfy a linear Weingarten
condition of type $\kappa _{1}=m\kappa _{2}+n$, where $m$ and $n$ are
constant and $\kappa _{1}$ and $\kappa _{2}$ denote the principal curvatures
at each point of the surface. We study the family of surfaces foliated by a uniparametric family of circles in parallel planes. We prove that the surface must be rotational or the surface is part of the family of Riemann examples of maximal surfaces ($m=-1$, $n=0$). Finally, we  consider the class of rotational surfaces for the case $n=0$, obtaining a first integration if the axis is timelike and spacelike and a complete description if the axis is lightlike.
\end{abstract}

\section{Introduction and results}\label{intro}

In Minkowski space there exists a family of maximal non-rotational surfaces foliated by circles in parallel planes
\cite{lls,lo1}. Maximal surfaces are spacelike surfaces whose mean curvature $H$ vanishes at every point of the surface. These surfaces play the same role of the classical Riemann examples of minimal surfaces in Euclidean space \cite{ni}. If one now assumes that the mean curvature $H$ is a non-zero constant, then any spacelike surface constructed by circles in parallel planes in Minkowski space is necessarily  a surface of revolution \cite{lo2} (in Euclidean space, the analogous result was proved by Nitsche in \cite{ni2}).

In this work we extend these results in a more general setting. We study spacelike
surfaces that satisfy a relation of type
\begin{equation}
\kappa _{1}=m\kappa _{2}+n\hspace*{.5cm}m\not=0,  \label{w1}
\end{equation}
where $\kappa _{1}$ and $\kappa _{2}$ are the principal curvatures of $M$
respectively, and $m$ and $n$ are constants. We say that $M$ is a \textit{%
linear Weingarten surface}. As  particular cases, this class of surfaces contains  the umbilical surfaces
when $(m,n)=(1,0)$ and  the surfaces with constant  mean curvature   if $m=-1$.
Linear Weingarten surfaces belong to a wider class of surfaces called Weingarten surfaces. A Weingarten surface is a surface that satisfies a smooth relation of type $W(\kappa_1,\kappa_2)=0$. Equation (\ref{w1}) describes the simplest case of function $W$, that is, that $W$ is linear in its variables. In Euclidean space, Weingarten surfaces are the focus of interest for many geometers, beginning in the fifties with works by Hopf, Chern and Hartman among others.  In Minkowski space, Weingarten surfaces have been studied in \cite{ag,bls,dk,ds1,kky}. The relation (\ref{w1}) in Euclidean space has been considered in
\cite{lo3,p1,p2}.

In the first part of the work, we consider linear Weingarten   surfaces foliated by a uniparametric family of circles in parallel planes. As we have previously said, for the particular case that  $(m,n)=(-1,0)$, that is, $H=0$,  there exist examples of non-rotational surfaces. We ask for the existence of examples of non-rotational surfaces foliated by circles in parallel planes that satisfy the  general relation (\ref{w1}). We prove:

\begin{theorem}
\label{t1} Let $M$ be a spacelike   surface in $\hbox{\bf E}_{1}^{3}$ foliated by circles  in parallel planes. If $%
M $ is a linear Weingarten surface, then $M$ is a surface of revolution or the surface is part of the family of Riemann examples of maximal surfaces.
\end{theorem}
Thus, if the surface is not rotational, then necessarily $H=0$ and the surface is one of the Riemann examples.

\begin{corollary}
Riemann examples of maximal surfaces are the only non-rotational spacelike  surfaces in $\e_1^3$ foliated  by circles in parallel planes that are linear Weingarten
surfaces.
\end{corollary}

Surfaces in different ambient spaces foliated by circles in parallel planes have studied in the literature \cite{ja,lo1,ni2,pa1,se}.

The second part of this article considers spacelike surfaces of revolution that are of linear Weingarten type. In such case equation (\ref{w1}) is an ordinary differential equation that determines the shape of the generating curve of the surface. Let $\hbox{\bf E}_{1}^{3}$ be the Minkowski three-dimensional space, that
is, the real vector space $\hbox{\bb R}^{3}$ endowed with the metric $\langle ,\rangle =(dx_{1})^{2}+(dx_{2})^{2}-(dx_{3})^{2}$,
where $(x_{1},x_{2},x_{3})$ denote the usual coordinates in $\hbox{\bb R}%
^{3} $. In
 $\e_1^3$ there are three types of rotational surfaces depending on the causal character of the
 axis of revolution. The equation (\ref{w1}) can not be integrate in all its generality, but we obtain a first integration of (\ref{w1}) when the axis is timelike or spacelike. If the axis is lightlike, we completely solve the equation.  Exactly,

\begin{theorem} Let $M$ be a spacelike rotational surface in $\e_1^3$ satisfying $\kappa_1=m\kappa_2$. After a rigid motion of the ambient space, the surface parametrizes as:
\begin{enumerate}
\item If the axis of revolution is timelike, then $X(u,v)=(u\cos(v),u\sin(v),z(u))$, with
$$z'(u)=\pm \frac{1}{\sqrt{1+c u^{-2/m}}},\ c>0.$$
\item If the axis of revolution is spacelike, then $X(u,v)=(u,z(u)\sinh(v),z(u)\cosh(v))$, with
$$z'(u)=\sqrt{1-c z(u)^{-2/m}}, \ c>0.$$
\item  If the axis of revolution is lightlike, then $X(u,v)=(-2uv,z(u)+u-uv^2,z(u)-u-uv^2)$, with
\begin{enumerate}
\item  $z(u)=c\log(u)$, $c>0$, if $m=2$.
\item   $z(u)=\frac{mc}{m-2}u^{\frac{m-2}{m}}$, $c>0$, if $m\not=2$.
\end{enumerate}
\end{enumerate}
\end{theorem}

\begin{remark} We point out  that the study  of this work can not carry for timelike surfaces.
A timelike surface in $\e_1^3$ is a surface whose induced metric is Lorentzian. In general,   the Weingarten endomorphism of a timelike surface is not diagonalizable and then  the relation (\ref{w1}) has not sense.
\end{remark}

\section{Preliminaries}


A vector $v\in\hbox{\bf E}_1^3$ is said spacelike if $\langle v,v\rangle>0$
or $v=0$, timelike if $\langle v,v\rangle<0$ and lightlike if $\langle
v,v\rangle=0$ and $v\not=0$. A submanifold $M\subset\e_1^3$  is said
spacelike, timelike or lightlike if the induced metric on $M$ is a
Riemannian metric (positive definite), a Lorentzian metric (a metric of
index $1$) or a degenerated metric, respectively. If $M$ is a straight-line $L=<v>$, this means that $v$ is spacelike, timelike or lightlike, respectively. If $M$ is a plane $P$, this is equivalent that
any orthogonal vector to $P$ is timelike, spacelike or lightlike
respectively. An immersion $x:M\rightarrow \hbox{\bf E}_{1}^{3}$ of a surface $M$
is called \textit{spacelike} if the induced metric $x^{\ast }\langle
,\rangle $ on $M$ is a Riemannian metric.

In Minkowski space $\hbox{\bf E}_1^3$ the
pseudohyperbolic surface plays the same role as a sphere in Euclidean space.  If $p_0\in\e_1^3$, the
pseudohyperbolic surface of radius $r>0$ and centered at $p_0$ is given by  $\hbox{\h}^{2,1}(r,p_0)=\{x\in\hbox{\bf E}_1^3;\langle x-p_0,x-p_0\rangle=-r^2\}$. From the Euclidean viewpoint, $\hbox{\bf H}^{2,1}(r)$ is the hyperboloid of
two sheets. So, if $O$ is the origin of coordinates,  $\hbox{\bf H}^{2,1}(r,O)$ satisfies the equation $x_1^2+x_2^2-x_3^2=-r^2$ which is obtained by rotating the
hyperbola $x_1^2-x_3^2=r^2$ in the plane $x_2=0$ with respect to the $x_3$%
-axis. This surface is spacelike with mean curvature $H=1/r$ and with Gauss
curvature $K=1/r^2$. Moreover it is an umbilical surface and so, it is a linear Weingarten surface.

We now describe the surfaces of revolution of $\e_1^3$. We study those rigid motions of the ambient space that leave
a straight-line $L$ pointwised fixed.    Depending on the axis $L$, there are three types of rotational motions. After an isometry of $\e_1^3$, the expressions of rotational motions with respect to the canonical basis  $\{e_1,e_2,e_3\}$  are as follows:
$$R_v:\left(\begin{array}{c}x_1\\ x_2\\ x_3\end{array}\right)\longmapsto
\left(\begin{array}{ccc}\cos{v}&\sin{v}&0\\-\sin{v}&\cos{v}&0\\0&0&1\end{array}\right)\left(\begin{array}{c}x_1\\ x_2\\ x_3\end{array}\right).$$
$$R_v: \left(\begin{array}{c}x_1\\ x_2\\  x_3\end{array}\right)\longmapsto
\left(\begin{array}{ccc}1&0&0\\0&\cosh{v}&\sinh{v}\\0&\sinh{v}&\cosh{v}\end{array}\right)\left(\begin{array}{c}x_1\\ x_2\\ x_3\end{array}\right).$$
$$R_v: \left(\begin{array}{c}x_1\\ x_2\\ x_3\end{array}\right)\longmapsto
\left(\begin{array}{ccc}1&-v&v\\ v&1-\frac{v^2}{2}&\frac{v^2}{2}\\ v&-\frac{v^2}{2}&1+\frac{v^2}{2}
\end{array}\right)\left(\begin{array}{c}x_1\\ x_2\\ x_3\end{array}\right).$$
A surface $M$ in $\e_1^3$ is a surface of revolution (or
rotational surface) if  $M$ is
invariant by some of the above three groups of rigid motions. In particular, there exists a planar curve $\alpha=\alpha(u)$ that generates the surface, that is, $M$ is the set of points given by
$\{R_v(\alpha(u));u\in I,v\in\r\}$. Because in this paper we are interested for spacelike surfaces, the generating curve $\alpha$ of the surface must be spacelike. We now describe the parametrizations of a spacelike rotational surface.
\begin{enumerate}
\item \emph{Case $L$ is a timelike axis}. Consider that  $L$ is the $x_3$-axis. If $p=(x,y,z)\not\in L$, then $\{R_v(p);v\in\r\}$ is an Euclidean circle of radius
$\sqrt{x^2+y^2}$ in the plane $x_3=z$. If $\alpha(u)=(u,0,z(u))$ is a planar curve in the plane $x_2=0$, then the surface of revolution generated by $\alpha$ writes as
\begin{equation}\label{rot1}
X(u,v)=(u\cos(v),u\sin(v),z(u)),\hspace*{1cm}z'^2<1, u\not=0.
\end{equation}

\item \emph{Case $L$ is a spacelike axis}.  Consider that $L$ is the $x_1$-axis. If $p=(x,y,z)$ does not belong to $L$, then $\{R_v(p);v\in\r\}$ is an Euclidean hyperbola in the
plane $x_1=x$ and with equation $x_2^2-x_3^2=y^2-z^2$.  If $\alpha(u)=(u,0,z(u))$ is a planar curve in the plane $x_2=0$, then the surface of revolution generated by $\alpha$ writes as
\begin{equation}\label{rot2}
X(u,v)=(u,z(u)\sinh(v),z(u)\cosh(v)),\hspace*{1cm}z'^2<1, u\not=0.
\end{equation}
\item \emph{Case $L$ is a lightlike axis}. Consider that $L$ is the  straight-line $v_1=<(0,1,1)>$. If $(x,y,z)$ does not belong to the plane $<e_1,v_1>$,    the orbit $\{R_v(p);v\in\r\}$ is the curve
 $$\beta(v)=(x-(y-z)v,xv+y-(y-z)\frac{v^2}{2},xv+z-(y-z)\frac{v^2}{2}).$$
 The curve $\beta$ lies in the plane $x_2-x_3=y-z$ and describes a parabola in this plane, namely,
$$\beta(v)=(x,y,z)+v(-(y-z) e_1+x v_1)-\frac{y-z}{2}v^2v_1.$$
Consider $\alpha(u)$ a planar curve in the plane $<(0,1,1),(0,1,-1)>$ given  as a graph on
the straight-line $<(0,1,-1)>$, that is, $\alpha(u)=(0,u+z(u),-u+z(u))$. The surface of revolution generated by $\alpha$ is
\begin{equation}\label{rot3}
X(u,v)=(-2 uv,z(u)+u -u v^2,z(u)-u-u v^2),\hspace*{1cm}z'>0, u\not=0.
\end{equation}
\end{enumerate}

Given a family of rotational motions, if  we look the orbit that describes a point of $p$ under the motions of this  family, we will obtain a planar curve that  plays the role of a circle in $\e_1^3$. From the above, and taking into account what happens in Euclidean ambient space, we give the definition of
a spacelike circle in Minkowski space $\e_1^3$.

\begin{definition}  An orbit of a point under one of the above three groups of rotations of $\e_1^3$ parametrized by a spacelike curve is called a spacelike circle.
\end{definition}
From now on, we say simply circle instead of spacelike circle. After a rigid motion of $\e_1^3$,
 there are three types of circles, which can viewed as Euclidean (horizontal) circles, (vertical) spacelike hyperbolas and spacelike parabolas.
Also one can prove that the definition of a circle    is equivalent to say a planar curve with constant curvature.

We end this section with some local computations for the curvatures of a spacelike surface, and that will be useful in next sections. Let $M$ be a spacelike surface in $\hbox{\bf E}_{1}^{3}$. The spacelike
condition is equivalent that any unit normal vector $\hbox{\bf N}$ to $M$ is
always timelike. Since any two timelike vectors in $\hbox{\bf E}_{1}^{3}$
can not be orthogonal, then we have $\langle \hbox{\bf N},(0,0,1)\rangle \neq 0$ on $M$.
This shows that $M$ is an orientable surface.

Let $x:M\rightarrow\e_1^3$ be a spacelike  immersion of a surface  $M$ and let  $N$ be a Gauss map.
Let $U, V$ be  vector fields to $M$ and we denote by  $\nabla^0$ and $\nabla$ the Levi-Civitta connections of  $\e_1^3$ and $M$ respectively.  The Gauss formula says
$\nabla_U^0 V=\nabla_U V+\mbox{II}(U,V),$
where $\mbox{II}$ is the second fundamental form of the immersion. The Weingarten endomorphism is $A_p:T_pM
\rightarrow T_p M$ defined as  $A_p(U)=-(\nabla_U^0 N)_p^\top=(-dN)_p(U)$.   We have then
 $\mbox{II}(U,V)=- \langle \mbox{II}(U,V),N\rangle N=- \langle AU,V\rangle N$. The mean curvature vector $\vec{H}$ is defined as
$\vec{H}=(1/2)\mbox{trace}(\mbox{II})$ and the Gauss curvature $K$ as the determinant of $\mbox{II}$ computed in both cases with respect to an orthonomal basis. The mean curvature $H$ is the function given by $\vec{H}=HN$, that is, $H=- \langle\vec{H},N\rangle$. If $\{e_1,e_2\}$ is an orthonormal vectors at each tangent plane,   then
$$\vec{H}=\frac12\mbox (\mbox{II}(e_1,e_1)+\mbox{II}(e_2,e_2))=- \frac12(\langle Ae_1,e_1\rangle+\langle Ae_2,e_2\rangle)N=- (\frac12 \mbox{trace}(A))N
$$
Then
\begin{equation}\label{hk}
H=-\frac{1}{2}\ \mbox{trace}(-d\hbox{\bf N})=-\frac{\kappa_1+\kappa_2}{2},\hspace*{1cm}K=-\mbox{det}\ (-d%
\hbox{\bf N})=-\kappa_1\kappa_2.\end{equation}
If we locally write the immersion as $X(u,v)$, with $(u,v)$ in
some planar domain, then
$$H=-\frac{1}{2}\ \frac{eG-2fF+gE}{EG-F^{2}},\hspace*{1cm}K=-\frac{e\ g-f^{2}}{%
EG-F^{2}},$$
where $\{E,F,G\}$ and $\{e,f,g\}$ are the coefficients of the first and
second fundamental forms respectively of the immersion according to the
orientation ${\bf N}=X_{u}\wedge X_{v}/|X_{u}\wedge X_{v}|$:
$$E=\langle X_u,X_u\rangle,\ F=\langle X_u,X_v\rangle,\ G=\langle X_v,x_v\rangle,$$
$$e=\langle N,X_{uu}\rangle,\ f=\langle N,X_{uv}\rangle,\ g=\langle N,X_{vv}\rangle,$$
where the subscripts denote the corresponding derivatives.
Denote $Q=EG-F^{2}=|X_{u}\wedge X_{v}|^{2}$. This
function is positive because the immersion is spacelike. From the
expressions of $H$ and $K$, we have
$$H_1:=-\Big(G[X_{u},X_{v},X_{uu}]-2F[X_{u},%
X_{v},X_{uv}]+E[X_{u},X_{v},%
X_{vv}]\Big)=2HQ^{3/2}  $$
$$K_1:=-\Big(\lbrack X_{u},X_{v},X_{uu}][X%
_{u},X_{v},X_{vv}]-[X_{u},X_{v},X%
_{uv}]^{2}\Big)=KQ^{2}$$
and $[,,]$ denotes the determinant of three vectors: $[v_{1},v_{2},v_{3}]=\mbox{det}(v_{1},v_{2},v_{3})$. The principal curvatures $\kappa _{1}$ and $\kappa _{2}$ are obtained by (\ref{hk}). So we have
$$\kappa _{1}=-H+\sqrt{H^{2}+K},\hspace*{.5cm} \kappa _{2}=-H-\sqrt{H^{2}+K}.$$
Then the condition (\ref{w1}) writes now as
$$(1-m)H_{1}+2nQ^{3/2}=(1+m)\sqrt{H_{1}^{2}+4QK_{1}},$$
and after some manipulations, and squaring twice the above equation, we obtain an expression of type
\begin{equation}\label{5}
(mH_{1}^{2}+(1+m)^{2}QK_{1}-n^{2}Q^{3})^{2}-n^{2}(1-m)^{2}H_{1}^{2}Q^{3}=0
\end{equation}


\section{Surfaces foliated by circles in parallel planes}


In this section we prove Theorem \ref{t1}. We consider a spacelike surface $M\subset \hbox{\bf E}_{1}^{3}$ parametrized
by circles in parallel planes. In the proof of the theorem we distinguish three cases according to the
causal character of the planes of the foliation.  Also we distinguish the case that the constant $n$ in (\ref{w1}) is or not zero. We also discard the case $m=1$, $n=0$, corresponding to the umbilical case (the surface is a pseudohyperbolic surface)  and the case  $m=-1$, where the surface has constant mean curvature:   it is known that the only spacelike   surfaces in $\e_1^3$ foliated by circles in parallel planes and with constant mean curvature are surfaces of revolution and the maximal Riemann examples (\cite{lo2}).  

\subsection{The planes are spacelike}

After a rigid motion in $\hbox{\bf E}_{1}^{3},$ we may assume the planes containing the circles of the foliation are
parallel to the plane $x_{3}=0$. According to the description given in (\ref{rot1}) in Preliminaries,  the   surface $M$ can be parametrized by
$$X(u,v)=(x(u),y(u),0)+(r(u)\cos {v},r(u)\sin {v},u),$$
where $x$, $y$ and $r$ are smooth functions in some interval $I\subset\r$, $r>0$. With this
parametrization, $M$ is a surface of revolution if and only if $x$ and $y$
are constant functions $x(u)=x_0$, $y(u)=y_0$ (the axis of revolution would be the straight line $x_1=x_0,x_2=y_0$).

The proof is by contradiction.  This means that the planar
curve $\alpha(u):=(x(u),y(u))$ is not constant, that is, $\alpha(u)$ is not a single point and we reparametrize $\alpha$
  by the arc-length, and thus $x'^2+y'^2=1$. Then  there exits a smooth function $\theta$ such that
$x'(u)=\cos\theta(u)$ and $y'(u)=\sin\theta(u)$. In fact, $\theta'(u)$ is the curvature $\kappa$ of $\alpha$. Equation (\ref{5}) is an expression of type
\begin{equation}
\sum_{j=0}^{12}A_{j}(u)\cos {(jv)}+B_{j}(u)\sin {(jv)}=0.  \label{6}
\end{equation}
Because the family of functions $\{\cos(jv),\sin(jv)\}$ are independent linear, the coefficients $A_j$ and $B_j$ must vanish in all its domain. We distinguish
two cases according to the value of $n.$

First, we assume $n\not=0$. The computation of $A_{12}$ and B$_{12}$ gives respectively:
$$A_{12}=\frac{n^4 r^{12}\cos(12\theta)}{2048},\hspace*{1cm}B_{12}=\frac{n^4r ^{12}\sin(12\theta)}{2048}.$$
Since $A_{12}=B_{12}=0$, we obtain $nr=0$: contradiction. Therefore, the proof for spacelike planes reduces to consider that  $n=0$ in the relation (\ref{w1}). Then the equation (\ref{5}) is $mH_{1}^{2}+(1+m)^{2}QK_{1}=0$ and the sum
 in (\ref{6}) is until $j=3$, with
\begin{eqnarray*}
A_{3} &=&-\frac{1}{4}(1+m)^{2}r(u)^5\kappa\sin(3\theta).\\
B_{3} &=&\frac{1}{4}(1+m)^{2}r(u)^{5}\kappa \cos(3\theta).
\end{eqnarray*}
 As $m+1\not=0$ and $r>0$, from $A_3=B_3=0$ we obtain $\kappa=0$. This means that the curve $\alpha$ is a straight-line and $\theta(u)=\theta_0$ is a constant function.   Taking into account this, we calculate the following coefficients $A_2$ and $B_2$:
$$A_2=\frac12\cos(2\theta_0)r^4(4mr'^2+(m+1)^2rr'').$$
$$B_2=\frac12\sin(2\theta_0)r^4(4mr'^2+(m+1)^2rr'').$$
Because $A_2=B_2=0$, we  have $4mr'^2+(m+1)^2rr''=0$.
The computations of $A_1$ and $B_1$ are
$$A_1=2\cos(\theta_0)r^4r'(2mr'^2+(1+m^2)rr'').$$
$$B_1=2\sin(\theta_0)r^4r'(2mr'^2+(1+m^2)rr'').$$
From $A_1=B_1=0$, we have $r'=0$ or $2mr'^2+(1+m^2) rr''=0$. If $r'=0$, then $r$ is a constant function. With this value of $r$, the computation of
$Q=EG-F^2$ gives $Q=-r^2\sin(\theta_0-v)^2<0$, which it is impossible. Thus, $2mr'^2+(1+m^2) rr''=0$. By combining with
$4mr'^2+(m+1)^2rr''=0$, and the fact that $m\not=0$, we obtain $r'^2+rr''=0$. From here, we have
$$A_2=\frac12(1-m)^2\cos(2\theta_0)r^4r'^2.$$
$$B_2=\frac12(1-m)^2\sin(2\theta_0)r^4r'^2.$$
Using that $m\not=1$,   we conclude that $r'=0$. Thus $r$ is a constant function, and we know then that this is  a contradiction.  This   shows the Theorem for    surfaces foliated by circles in spacelike parallel planes.

\subsection{The planes are timelike}

We consider a spacelike surface foliated by circle in parallel timelike planes. After a rigid motion of the ambient space $\e_1^3$, we assume that these planes are parallel to the plane $x_1=0$. By (\ref{rot2}), the surface locally parametrizes as
$$X(u,v)=(0,y(u),z(u))+(u,r(u)\sinh(v),r(u)\cosh(v)),$$
with $y, z$ and $r$ smooth functions, with $r>0$. The surface is a surface of revolution if the curve $\alpha(u)=
(0,y(u),z(u))$ is a constant point. Again, the proof is by contradiction, and we assume that $\alpha$ is not constant. We parametrize $\alpha$ so that $y'^2-z'^2=1$. This means that $y'(u)=\cosh\theta(u)$ and $z'(u)=\sinh\theta(u)$, for some function $\theta$. Moreover, $\theta'=\kappa$ is the curvature of $\alpha$. Now Equation (\ref{5}) writes as
\begin{equation}
\sum_{j=0}^{12}A_{j}(u)\cosh {(jv)}+B_{j}(u)\sinh {(jv)}=0.  \label{62}
\end{equation}
The functions  $\{\cosh(jv),\sinh(jv)\}$ are independent linear and thus the coefficients $A_j$ and $B_j$ must vanish in all its domain. We calculate  $A_{12}$ and $B_{12}$, obtaining
$$A_{12}=\frac{n^4 r^{12}\cosh(12\theta)}{2048},\hspace*{1cm}B_{12}=-\frac{n^4r ^{12}\sinh(12\theta)}{2048}.$$
Since $A_{12}=B_{12}=0$ and $n\not=0$, we obtain a contradiction.

Therefore, we suppose that $n=0$. Then  the sum
 in (\ref{62}) is until $j=3$ again, with
\begin{eqnarray*}
A_{3} &=&\frac{1}{4}(1+m)^{2}r(u)^5\kappa\cosh(3\theta).\\
B_{3} &=&-\frac{1}{4}(1+m)^{2}r(u)^{5}\kappa \sinh(3\theta).
\end{eqnarray*}
Since $m+1\not=0$ and $r>0$, from $A_3=B_3=0$ we obtain $\kappa=0$. As in the previous case, $\alpha$ is a straight-line and $\theta(u)=\theta_0$ is a constant function.   The coefficients  $A_2$ and $B_2$ are
$$A_2=\frac12\cosh(2\theta_0)r^4(4mr'^2+(m+1)^2rr'').$$
$$B_2=-\frac12\sinh(2\theta_0)r^4(4mr'^2+(m+1)^2rr'').$$
The coefficients $A_1$ and $B_1$ are
$$A_1=2\sinh(\theta_0)r^4r'(-4m+2mr'^2+(1+m^2)rr'').$$
$$B_1=2\cosh(\theta_0)r^4r'(-4m+2mr'^2+(1+m^2)rr'').$$
From $A_1=B_1=0$, we have $r'=0$ or $-4m+2mr'^2+(1+m^2) rr''=0$. Assume $r'=0$. Then the coefficient $A_0$ is
$A_0=4mr^4$. As $A_0=0$ and $r\not=0$, we get a contradiction. Suppose now $-4m+2mr'^2+(1+m^2) rr''=0$. By combining with
$4mr'^2+(m+1)^2rr''=0$, we have $r'^2+2+rr''=0$. We use it in $A_1$  obtaining $(m-1)^2r'^2+2(m+1)^2=0$.
In particular, $m=-1$, which it was discarded at the beginning of this section. This contradiction shows the result.
\subsection{The planes are lightlike}

After a motion in $\e_1^3$ we may assume the foliating planes
are parallel to $x_{2}-x_{3}=0$. Recall that now a circle is a parabola whose axis is parallel to the vector
$(0,1,1)$. If the surface is rotational, the vertex of the parabolas $v\longmapsto X(u,v)$ in (\ref{rot3}) is given by $\alpha(u)$. However, in the case that the surface is not rotational,  we allow that this vertex belongs the plane $x_2-x_3$. Then the surface is parametrized by
$$X(u,v)=(a(u),0,0)+(-2uv,b(u)+u-uv^2,b(u)-u-uv^2),$$
where  $a$ and $b$  are smooth functions. The surface $M$ will be rotational if the function $a(u)$ is constant.
The equation (\ref{5}) writes now as
\begin{equation}\label{light}
\sum_{j=0}^{6}A_{j}(u)v^{n}=0,
\end{equation}
for smooth functions $A_{j}$. If we look (\ref{light}) as polynomial on $v$, then all coefficients $A_{j}(u)$
vanish. Assume first that $n\not=0$.  Then $A_{6}=0$ means $n^4u^{12} a'^6=0$ and so, $a'=0$, that is,
$a$ is a constant function.

Suppose now $n=0$ in the relation (\ref{w1}). Now we use (\ref{5}) and the sum  (\ref{light}) is now until
$j=2$. We have
$$A_2=256u^4(2ma'+ua'')(2a'+mua'')=0.$$
We assume that $2ma'+ua''=0$ (the reasoning with $2a'+mua''=0$ is similar). If $a$ is not a constant function, then
$a'(u)=cu^{-2m}$, with $c>0$. Putting now in $A_1$ we have
\begin{equation}\label{p1}
A_1=512c u^{4-6m}(m+1)\Big(-mc^2+(m-1)u^{4m}(2mb'+ub'')\Big).
\end{equation}
As $m\not=-1$, we obtain a value of $c$, which substituted into $A_0$ gives
$$(m+1)^2u^{4}(2mb'+ub'')^2=0,$$
and so, $2mb'+ub''=0$. Thus (\ref{p1}) implies now $mc^2=0$: contradiction. As conclusion,  $a=a(u)$ is a constant function.

\section{Rotational surfaces that satisfy $\kappa_1=m\kappa_2$}

In this section we study rotational spacelike surfaces that satisfy the relation (\ref{w1}). Due to the complexity in the general case, and the impossibility to obtain a complete integration of this equation, we focus in (\ref{w1}) when $n=0$. Then the Weingarten relation writes as
\begin{equation}\label{n0}
mH_1^2+(1+m)^2QK_1=0.\end{equation}
We distinguish the classification according the causal character of the axis of revolution.

\subsection{The axis is timelike}

Assume that the axis of revolution is timelike. According to (\ref{rot1}), we parametrize the surface
as $X(u,v)=(u \cos(v),u\sin(v),z(u))$, with  $z'^2<1$. A straightforward computation implies that identity (\ref{n0}) writes as
\begin{equation}\label{timelike}
-z'(1-z'^2)+m u z''=0\hspace*{.5cm}\mbox{or}\hspace*{.5cm}-mz'(1-z'^2)+uz''=0.
\end{equation}
Both equations are equivalent since they describe the situation $\kappa_1=m\kappa_2$ and $\kappa_2=m\kappa_1$, respectively. If $z'(u_0)=0$ at some point $u_0>0$, then the function $z(u)=z(u_0)$ is the solution, that is, the surface is a horizontal spacelike plane. Thus we suppose that $z'\not=0$ at some point and we consider the first equation in (\ref{timelike}). Let $\varphi=z'$. Then (\ref{timelike}) is
$$-\varphi(1-\varphi^2)+m u \varphi'=0.$$
A first integration leads to
$$\frac{\varphi}{\sqrt{1-\varphi^2}}=\mu u^{1/m},\ \mu>0.$$
Thus
$$z'=\pm \frac{1}{\sqrt{1+c u^{-2/m}}},$$
where $c$ is a positive integration constant.   The solution of this differential equation are given in terms of Gauss hypergeometric functions (see \cite[Ch. 15]{ab}).
 It is known that hypergeometric functions with special arguments reduce to elementary functions. This is the 
 case if  $m=1$ and $m=-1$. In fact,
if $m=1$, the solution is $z(u)=\lambda\pm \sqrt{u^2+c}$ and the surface is   the pseudohyperbolic surface $\h^{2,1}(\sqrt{c},p_0)$, with $p_0=(0,0,\lambda)$.  If $m=-1$, then $z(u)=\frac{1}{\sqrt{c}}\mbox{arc}\sinh{( \sqrt{c}u)+\lambda}$, $c>0$, $\lambda\in\r$, that is, the surface is the Lorentzian catenoid of first kind \cite{ko}. For other cases, it is possible to obtain an explicit solution of $z$. For example, for $m=\pm 2$ and letting $c=1$, we have:
\begin{enumerate}
\item $m=2$, $z(u)=\sqrt{u(1+u)}-\mbox{arc}\sinh(\sqrt{u})+\lambda$, $\lambda\in\r$.
\item $m=-2$, $z(u)=2\sqrt{1+u}+\lambda$, $\lambda\in\r$.
\end{enumerate}

\subsection{The axis is spacelike}

We  assume  from (\ref{rot2}) that the surface parametrizes as $X(u,v)=(u,z(u)\sinh(v),z(u)\cosh(v))$,
 with  $z'^2<1$.
Equation (\ref{n0}) writes as
\begin{equation}\label{s1}
-1+z'^2+m z z''=0\hspace*{.5cm}\mbox{or}\hspace*{.5cm}-m(1-z'^2)+z'z''=0.
\end{equation}
Letting $\varphi=z'$ and $\zeta=z$ as new dependent and independent variables, respectively, we only take the first equation in (\ref{s1}), since the second one is the analogous for $\kappa_2=m\kappa_1$. Then it transforms into
$$-1+\varphi^2+m\zeta \varphi\varphi'=0.$$
 A first integration leads to:
$$\frac{1}{\sqrt{1-z'^2}}=\mu z^{1/m}$$
for some positive constant $\mu$, or equivalently,
$$z'=\sqrt{1-c z^{-2/m}},\ c>0.$$
Again, the solutions are given in terms of hypergeometric functions. Let us see some exact solutions for special choices of the parameter $m$. For $m=1$, $z(u)=\pm\sqrt{(u+\lambda)^2+c}$, $\lambda\in\r$
 and for $m=-1$, $z(u)=\frac{1}{\sqrt{c}}\sin(\sqrt{c}u)+\lambda$, $\lambda\in\r$. The first example, it is the pseudohyperbolic surface $\h^{2,1}(\sqrt{c},p_0)$, $p_0=(-\lambda,0,0)$
and the second one  is a maximal surface called the catenoid of second kind \cite{ko}. For some other values of $m$, we can obtain explicit integrations. For example, taking $c=1$ and  $m=-2$, we obtain  $z(u)=\frac{4-(u+\lambda)^2}{4}$, $\lambda\in\r$.

\subsection{The axis is lightlike}
After a rigid motion of $\e_1^3$, we parametrize the surface as $X(u,v)=(-2 uv,z(u)+u -u v^2,z(u)-u-u v^2)$.  Equation (\ref{n0}) writes as
\begin{equation}\label{lightlike}
 2z'+muz''=0\hspace*{0.5cm}\mbox{or}\hspace*{0.5cm}2mz'+u z''=0.
 \end{equation}
Again, and as in the previous cases, we only consider the first equation. A first integration  leads $z'(u)=cu^{-2/m}$, with $c>0$. The complete integration of (\ref{lightlike}) gives
\begin{enumerate}
\item If $m=2$, then $z(u)=c\log(u)+\lambda$, $\lambda\in\r$.
\item If $m\not=2$, then $z(u)=\frac{mc}{m-2}u^{\frac{m-2}{m}}+\lambda$, $\lambda\in\r$.
\end{enumerate}
For the case that $m=1$, we have $z(u)=-c/u+\lambda$ and the coordinates of the surface satisfies $x_1^2+(x_2-\lambda)^2-(x_3-\lambda)^2=-4c$. Thus the surface is the pseudohyperbolic surface $\h^{2,1}(2\sqrt{c},p_0)$, with $p_0=(0,\lambda,\lambda)$. If $m=-1$, then $z(u)=cu^{3}/3+\lambda$ and the surface is the Enneper surface of second kind \cite{ko}.

\begin{figure}[hbtp]
\begin{center}\includegraphics[width=6cm]{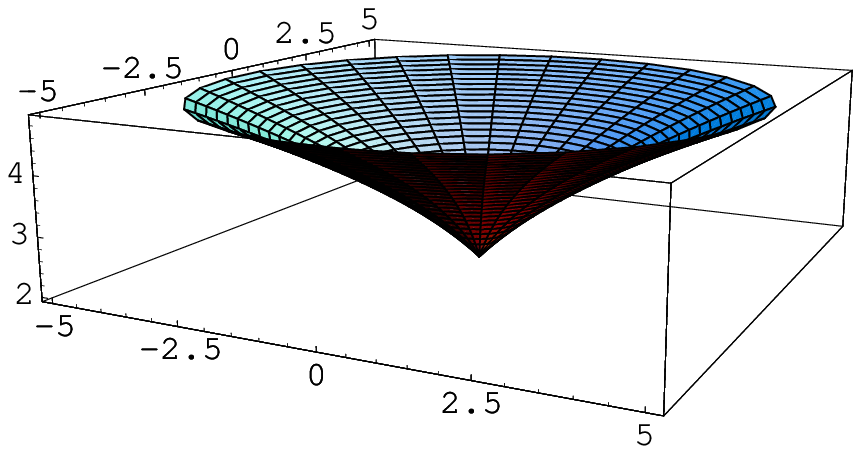}\hspace*{.5cm}\includegraphics[width=3cm]{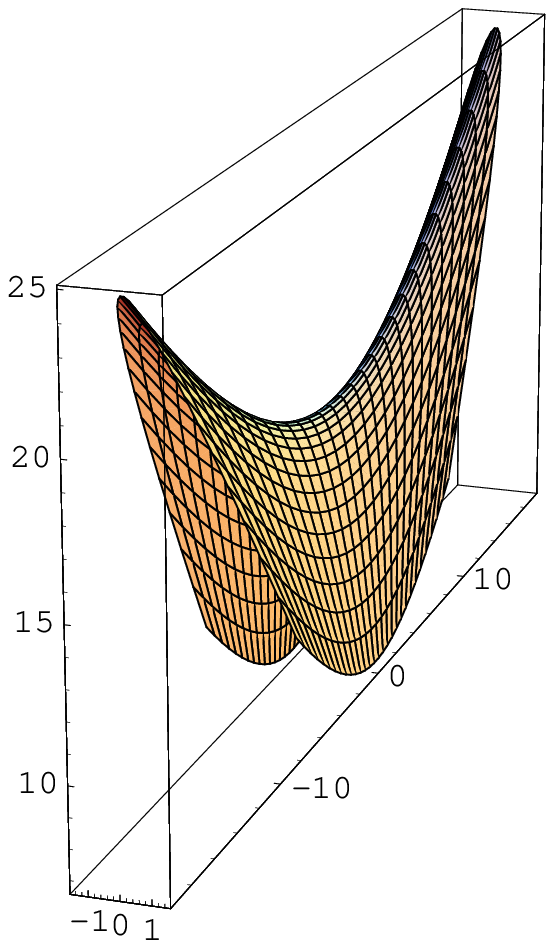}\hspace*{.5cm}
\includegraphics[width=5cm]{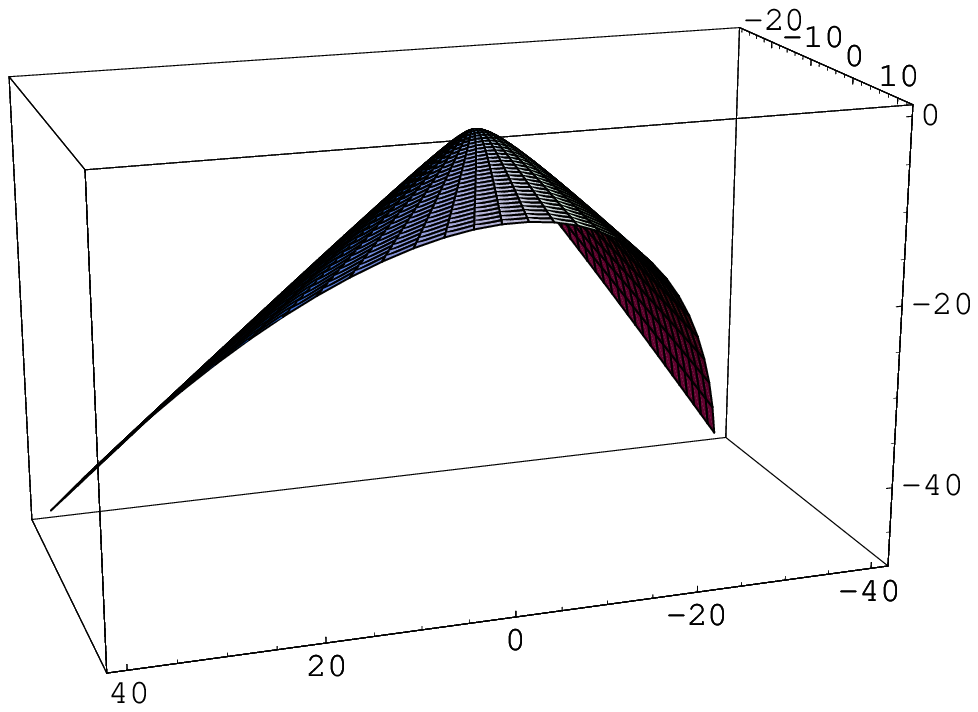}
\end{center}
\caption{Rotational spacelike surfaces that satisfy the relation (\ref{n0}) for $m=-2$: (left) timelike axis; (center) spacelike axis; (right) lightlike axis.}
\end{figure}

\end{document}